\documentclass[12pt]{amsart}
\usepackage{amssymb}
\newtheorem{theorem}{Theorem}[section]
\newtheorem{lemma}[theorem]{Lemma}

\newtheorem{definition}[theorem]{Definition}

\newtheorem{corollary}[theorem]{Corollary}
\newcommand{\M}{\mbox{$\mathbb{M}$}}
\newcommand{\C}{\mbox{$\mathbb{C}$}}
\newcommand{\N}{\mbox{$\mathbb{N}$}}
\newcommand{\K}{\mbox{${\mathcal K}$}}
\newcommand{\B}{\mbox{$\mathbb{B}$}}
\newcommand{\Hi}{\mbox{${\mathcal H}$}}

\newcommand{\Sy}{\mbox{${\mathcal S}$}}

\newcommand{\Q}{\mbox{${\mathcal{QC}}$}}

\newcommand{\A}{\mbox{${\mathcal A}$}}
\newcommand{\QM}{\mbox{${\mathcal{QM}}$}}
\newcommand{\Bra}{\mbox{${\mathcal B}$}}
\newcommand{\cl}[1]{\mathcal{#1}}

\begin{document}
\title[Quasi-multipliers and algebrizations of an operator space]
{Quasi-multipliers and algebrizations of an operator space}
\author{Masayoshi Kaneda*}
\address{Masayoshi KANEDA: Department of Mathematics\\
University of California, Irvine, CA 92697-3875 U.S.A.}
\email{mkaneda@math.uci.edu, http://www.math.uci.edu/$\thicksim$mkaneda/}
\date{\today}
\thanks{THE MAJOR RESULTS IN THIS PAPER WERE PRESENTED AT THE 23RD ANNUAL GREAT
PLAINS OPERATOR THEORY SYMPOSIUM, UNIVERSITY OF ILLINOIS, URBANA-CHAMPAIGN, May
29, 2003}
\thanks{* This research was supported by a grant from the National Science
Foundation}
\begin{abstract}
Let $X$ be an operator space, let $\varphi$ be a product on $X$, and let $(X,
\varphi)$ denote the algebra that one obtains. We give necessary and sufficient
conditions on the bilinear mapping $\varphi$ for the algebra $(X, \varphi)$ to
have a completely isometric representation as an algebra of operators on some
Hilbert space. In particular, we give an elegant geometrical characterization
of such products by using the Haagerup tensor product. Our result makes no
assumptions about identities or approximate identities. Our proof is
independent of the earlier result of Blecher-Ruan-Sinclair (\cite{BRS}) that
solved the case when the algebra has an identity of norm one, and our result is
used to give a simple direct proof of this earlier result. We also develop
further the connections between quasi-multipliers of operator spaces, and shows
that the quasi-multipliers of operator spaces defined in \cite{KP} coincide
with their $C^*$-algebraic counterparts.


\end{abstract}
\maketitle
\section{Introduction.}\label{section: intro}
One of the most interesting questions in the operator space theory was: what
are the possible operator algebra products which a given operator space can be
equipped with? I was investigating many types of multipliers of operator
spaces for their own interests. Meantime, V. I. Paulsen defined
quasi-multipliers of operator spaces (\cite{KP} Definition 2.2), and suggested
to me to study them. Then, accidentally, I found that the quasi-multipliers
happened to answer the question above. That is, the possible operator algebra
products which a given operator space can be equipped with are precisely the
bilinear mappings implemented by the contractive quasi-multipliers of the
operator space (\cite{KP} Theorem 2.6). In this paper, we give a striking
geometrical characterization of operator algebra products (Theorem~\ref{main}).

In Section~\ref{section: two-sided}, we study quasi-multipliers in a special
case in which an operator space is an operator algebra with a two-sided
contractive approximate identity (we will abbreviate as ``c.a.i.''). In this
case, the quasi-multiplier space is quite manageable like $C^*$-algebra case,
and equivalent to other definitions using representation on a Hilbert space, or
considering in the second dual. It is also equivalent to the set of
quasi-centralizers, as a result, we obtain that the definition of
quasi-multipliers of operator spaces coincide with the existing ones in the
$C^*$-algebra case which were defined by L.~G.~Brown (\cite{Br}).

In Section~\ref{section: algebrizations}, we present the main result (Theorem
\ref{main}) of this paper. There, we give a beautiful geometrical
characterization of operator algebra products under no assumption about
identities or approximate identities. That is, the possible operator
algebra products which a given operator space can be equipped with are
completely determined by matrix norm structure of the operator space by using
the Haagerup tensor product (Theorem \ref{main}). This can be considered as
the quasi-version of the Blecher-Effros-Zarikian theorem ($\tau$-trick) in
which they characterized left multiplier mappings in terms of the matrix norms
(\cite{BEZ}). As a simple corollary, we obtain a generalized version
of Blecher-Ruan-Sinclair theorem (\cite{BRS}).

The reader who hurries for the main result may skip Section~\ref{section:
two-sided} and directly move on to Section~\ref{section: algebrizations} after
reading Section~\ref{section: pre} for a back ground if necessary.
\section{Preliminaries.}\label{section: pre}
We begin by recalling a construction of an injective envelope of an operator
space. See, e.g., \cite{BP}, \cite{P} Chapter 15 for more details. Let
$X\subset\B(\K, \Hi)$ be an operator
space, and consider the Paulsen operator system$$\Sy_X:=\left[\begin{matrix}
\C1_{\Hi}&X\\X^*&\C1_{\K}\\\end{matrix}\right]\subset\B(\Hi\oplus\K).$$One then
takes a minimal (with respect to a certain ordering) completely positive
$\Sy_X$-projection $\Phi$ on $\B(\Hi\oplus\K)$ whose image Im$\Phi$ turned out
to be an injective envelope $I(\Sy_X)$ of $\Sy_X$. By a well-known result of
M.-D.~Choi and E.~G.~Effros (\cite{CE}), Im$\Phi$ is a unital C*-algebra with
the product $\odot$ defined by $\xi\odot\eta:=\Phi(\xi\eta)$ for $\xi,
\eta\in$Im$\Phi$ and other algebraic operations and norm are the original ones
in $\B(\Hi\oplus\K)$. One may write$$\text{Im}\Phi=I(\Sy_X)=\left[
\begin{matrix}I_{11}(X)&I(X)\\I(X)^*&I_{22}(X)\end{matrix}\right]\subset\B(\Hi
\oplus\K),$$where $I(X)$ is an injective envelope of $X$, and $I_{11}(X)$ and
$I_{22}(X)$ are injective unital $C^*$-algebras.

By well-known trick one may decompose$$\Phi=\left[\begin{matrix}\psi_1&\phi\\
\phi^*&\psi_2\end{matrix}\right].$$

The new product $\odot$ induces new products $\bullet$ between elements of
$I_{11}(X)$, $I_{22}(X)$, $I(X)$ and $I(X)^*$. For example, $x\bullet a=
\phi(xa)$ for $x\in I(X)$, $a\in I_{22}(X)$. Note that the associativity of
$\bullet$ is guaranteed by that of $\odot$.

We call the embedding $i: X\hookrightarrow\left[\begin{matrix}O&X\\O&O
\end{matrix}\right]\subset\Sy_X;\;x\to\left[\begin{matrix}0&x\\0&0\end{matrix}
\right]$ the {\em \v{S}ilov embedding}.

\vspace{5mm}
Now we recall the definitions of quasi-multipliers for operator spaces, and
also define some related notions.
\begin{definition}\label{qm}
\begin{itemize}
\item[(1)] (\cite{KP} Section 2) Let $X$ be an operator space, and let $\pi$ be
a complete isometry from $X$ into an operator algebra $\A$. Then
{\bf $(\A, \pi)$-relative quasi-multiplier space} of $X$ is the set
$$\QM^{\pi}(X):=\{a\in\A;\;\pi(X)a\pi(X)\subset\pi(X)\}.$$
\item[(2)] Let $\QM^{\pi}(X)$ and $\QM^{\pi'}(X)$ be, respectively,
$(\A, \pi)$-relative and $(\A', \pi')$-relative quasi-multiplier spaces for
$X$. Then a linear mapping $\sigma: \QM^{\pi}(X)\to\QM^{\pi'}(X)$ is a {\bf
quasi-homomorphism}\footnote{This is different from a {\em left (right)
quasihomomorphism} defined in \cite{KP} Definition 5.1.} if $\pi^{-1}(\pi(x_1)y
\pi(x_2))=\pi'^{-1}(\pi'(x_1)\sigma(y)\pi'(x_2)),\;\;\forall x_1, x_2\in X,\;y
\in\QM^{\pi}(X)$. Furthermore, if $\sigma$ is one-to-one and onto, then we call
$\sigma$ a {\bf quasi-isomorphism}.
\item[(3)](\cite{KP} Definition 2.6) The {\bf quasi-multiplier space} for $X$
is the set $$\QM(X):=\{z\in I(X)^*;\;\;X\bullet z\bullet X\subset X\}.$$ We
call an element of $\QM(X)$ a {\bf quasi-multiplier} of $X$.
\end{itemize}
\end{definition}
Note that $\QM(X)$ is a subspace of $\QM^i(X)$ under the identification\\
$\QM(X)=\left[\begin{matrix}O & O\\\QM(X) & O\end{matrix}\right]$, where $i$ is
the \v{S}ilov embedding $X\to\left[\begin{matrix}O & X\\O &
O\end{matrix}\right]$ defined above.

\vspace{5mm}
The following theorem shows the universal property of quasi-multipliers.
\begin{theorem}\label{situation} (\cite{KP} Theorem 2.3) Let $X$ be an operator
 space and $\A$ be an operator algebra, and suppose that $\pi: X\to\A$ is a
complete isometry. Then there exists a unique completely contractive
quasi-homomorphism $\sigma: \QM^{\pi}(X)\to\QM(X)$, i.e., $\pi(x_1)y\pi(x_2)=
\pi(x_1\bullet \sigma(y)\bullet x_2),\;\;\forall x_1, x_2\in X,\;y\in\QM^{\pi}
(X)$, where $X$ is regarded as a subset of $I(\Sy_X)$. In particular, if $X$
itself is an operator algebra with product $\cdot$, then there exists a unique
$z\in\QM(X)$ such that $x_1\cdot x_2=x_1\bullet z\bullet x_2,\;\;\forall x_1,
x_2\in X$.
\end{theorem}

\vspace{5mm}
In this paper, we discuss bilinear mappings on operator spaces. We would like
to make sure of our terminology. Let $\varphi$ be a bilinear mapping on an
operator space $X$, and let $\tilde{\varphi}: X\otimes_hX\to X$ be the linear
mapping corresponding to $\varphi$, where $X\otimes_hX$ is the Haagerup tensor
product which plays the central role in Section \ref{section: algebrizations}.
For the Haagerup tensor product, see \cite{P}, \cite{ER}, \cite{Pi}. We define
the {\bf completely bounded norm} of $\varphi$ by
$\|\varphi\|_{cb}:=\|\tilde{\varphi}\|_{cb}$, and we say that $\varphi$ is
{\bf completely bounded} (respectively, {\bf completely contractive}) if
$\|\varphi\|_{cb}<\infty$ (respectively, $\|\varphi\|_{cb}\le1$). Note that the
term ``completely bounded'' for a bilinear mapping defined here is in the sense
of Christensen-Sinclair \cite{CS}, and it is called {\em multiplicatively
bounded} in \cite{ER}.

Throughout this paper, a {\bf product} means an {\em associative} bilinear
mapping.
\section{Quasi-multipliers of an operator algebra with a two-sided c.a.i.}
\label{section: two-sided}
In this section, we study some equivalent notions of quasi-multipliers of an
operator algebra with a two-sided c.a.i. in a similar manner to \cite{B},
\cite{BK}. We have
already studied the left multipliers (\cite{BK}) and the right multipliers
(\cite{B}) of an operator algebra with a ``right'' c.a.i.. So the reader may
think that we should study here the quasi-multipliers of an operator algebra
with a right c.a.i.. However, this seems very messy and not worth writing here,
so we study only the ``two-sided c.a.i.'' case. This is very satisfactory, and
as a result, we show that ``quasi-multipliers'' coincide with the existing ones
in the C*-algebra case (\cite{Br}, also see \cite{Pe} \S3.12).
\begin{definition}
Let $\A$ be an operator algebra with a two-sided c.a.i.. Then
the {\bf quasi-centralizer space} of $\A$ is the set
$$\Q(\A):=\{\varphi: \A\times\A\to\A;\;\;\varphi(\cdot, a)\in{_{\A}CB(\A)},\;\;
\varphi(a, \cdot)\in CB_{\A}(\A)\;\;\forall a\in\A\}.$$ We
call an element of $\Q(\A)$ a {\bf quasi-centralizer}. In Theorem \ref{qmti},
we prove that elements of $\Q(\A)$ are completely bounded in the sense of
Christensen-Sinclair which is explained at the end of the previous section.
Hence, for each $n\in\N$ and each $(\varphi_{i, j})\in\M_n(\Q(\A))$, there is a
corresponding matrix of linear mappings $(\tilde{\varphi}_{i, j})\in\M_n(CB(\A
\otimes_h\A, \A))$, and we define a matrix norm of $\Q(\A)$ by $\|(\varphi_{i,
j})\|_n:=\|(\tilde{\varphi}_{i, j})\|_n$.

We now define a {\bf quasi-multiplier extension} of $\A$ to be a
pair $(X, \pi)$ consisting of an operator space $X$ which is a subspace of some
operator algebra $\Bra$ and a completely isometric homomorphism $\pi: \A\to X$
such that $\pi(\A)X\pi(\A)\subset\pi(\A)$, where the product is taken in
$\Bra$. We say that $(X,
\pi)$ is an {\bf essential} quasi-multiplier extension of $\A$ if in addition
the canonical completely contractive mapping $X\to\Q(\A)$ is one-to-one. For
two quasi-multiplier extension $(X, \pi)$ and $(X', \pi')$ of $\A$, we write
$(X, \pi)\leq(X', \pi')$ if there exists a completely contractive homomorphism
$\theta: X\to X'$ such that $\theta\circ\pi=\pi'$. 
We say that two quasi-multiplier extensions $(X, \pi)$ and
$(X', \pi')$ are {\bf $\A$-equivalent} if there exists a completely
isometric quasi-isomorphism $\theta: X\to X'$ with $\theta\circ\pi=\pi'$. This
is an equivalence relation, and ``$\le$'' induces a well-defined ordering
on the equivalence classes.

\end{definition}
It follows that if there exists a maximum essential quasi-multiplier extension
of $\A$, then it is unique up to $\A$-equivalence. Also if two quasi-multiplier
extensions are $\A$-equivalent, and if one is essential, then so is the other.

\vspace{5mm}
Now let us recall the following facts from \cite{BK}.
\begin{lemma}\label{bk}(\cite{BK} Theorem 2.3) Let $\A$ be an operator algebra
with a right c.a.i.. Then there exists a $v^*\in Ball(\QM(\A))$ such that
\begin{enumerate}
\item $v\bullet v^*$ is an orthogonal projection in $I_{11}(\A)$,
\item $v^*\bullet v$ is the identity of $I_{22}(\A)$,
\item $a\bullet v^*\bullet b=ab,\;\forall a,b\in\A$, and hence $\hat{\psi}: I(
\A)\to I_{11}(\A)$ defined by $\hat{\psi}(a):=a\bullet v^*,\;\forall a\in I(\A)
$ is a complete isometry that restricts to a homomorphism on $\A$.
\end{enumerate}
\end{lemma}
The following corollary immediately follows from the lemma above.
\begin{corollary}\label{corollary: bk}If $\A$ has a two-sided contractive
approximate identity, then $v\bullet v^*$ is the identity of $I_{11}(\A)$.
\end{corollary}
\begin{proof}By symmetry, there exists a $w^*Ball(\QM(\A))$ such that
$w\bullet w^*$ is the identity of $I_{11}(\A)$. But $v^*=w^*$ by uniqueness of
a quasi-multiplier (Theorem \ref{situation}).
\end{proof}
The following lemma tells us that if an operator algebra $\A$ has a two-sided
c.a.i., then the quasi-multiplier space $\QM(\A)$ can be replaced by a better
one, $\QM^{\psi}(\A)$, which contains a copy of $\A$ preserving the product.
\begin{lemma}\label{improvedqm} Let $\A$ be an operator algebra with a
two-sided c.a.i., and we regard $\A$ as $\A=\left[\begin{matrix}O&\A\\O&O
\end{matrix}\right]\subset I(\Sy_{\A})$. Let $\tilde{\psi}: \A\to I_{11}(\A);\;
\tilde{\psi}:=\hat{\psi}|_{\A}$ and $\rho: \QM(\A)\to I_{11}(\A);\;\rho(z):=v
\bullet z$, where $\hat{\psi}$ and $v$ are as in Lemma \ref{bk}. Then $\rho$ is
a completely isometric quasi-isomorphism from $\QM(\A)$ onto $\QM^{\tilde{\psi}
}(\A)$, and $\tilde{\psi}(\A)\subset\QM^{\tilde{\psi}}(\A)\;(=\rho(\QM(\A)))$.
\end{lemma}
\begin{proof} For $a_1, a_2\in\A,\;z\in\QM(\A)$,\;
$\tilde{\psi}(a_1)\bullet\rho(z)\bullet\tilde{\psi}(a_2)=a_1\bullet v^*\bullet
v\bullet z\bullet a_2\bullet v^*=a_1\bullet z\bullet a_2\bullet
v^*=\tilde{\psi}(a_1\bullet z\bullet a_2)$.
This shows that $\rho$ is a quasi-homomorphism and
$\rho(\QM(\A))\subset\QM^{\tilde{\psi}}(\A)$. On the contrary, let $x\in
I_{11}(\A)$ be such that $\tilde{\psi}(a_1)\bullet
x\bullet\tilde{\psi}(a_2)=\tilde{\psi}(a_3)$ for some $a_3\in\A$. Then
$a_1\bullet v^*\bullet x\bullet a_2\bullet v^*=a_3\bullet v^*$. By multiplying
the both sides by $v^*$ on the right, we obtain $a_1\bullet v^*\bullet x\bullet
a_2=a_3$. This implies that $v^*\bullet x\in\QM(\A)$. Hence, $x=v\bullet v^*
\bullet x\in\rho(\QM(\A))$ since $v\bullet v^*$ is the identity of $I_{11}(\A)$
(Corollary \ref{corollary: bk}). That $\rho$ is a complete isometry easily
follows from the facts that $\rho$ is a left multiplication by $v$ and that
$v^*\bullet v$ is the identity of $I_{22}(\A)$. That $\tilde{\psi}(\A)\subset
\QM^{\tilde{\psi}}(\A)$ follows from the fact that $\tilde{\psi}$ is a
homomorphism.
\end{proof}
Let $\A$ be an operator algebra with a two-sided c.a.i.. We consider the
following three notions:
\begin{itemize}
\item[(I)] $\QM^{**}(\A):=\{z\in \A^{**};\;
\cl{\hat{A}}z\cl{\hat{A}}\subset\cl{\hat{A}}\}$,
\item[(II)] $\QM^{\pi}(\A):=\{T\in\B(\Hi); \pi(\A)T\pi(\A)\subset\pi(\A)\}$,
where $\pi: \A\to\B(\Hi)$ is a completely isometric nondegenerate
representation,
\item[(III)] $\QM^{\tilde{\psi}}(\A)$,
\end{itemize}
where $\cl{\hat{A}}$ is the canonical image of $\A$ in $\A^{**}$. Also
hereafter, we denote by $\hat{a}\in\cl{\hat{A}}$ the canonical image of
$a\in\A$.
\begin{theorem}\label{qmti} Let $\A$ be an operator algebra with a two-sided
c.a.i.. Then (I), (II), and (III) as above are quasi-multiplier extensions of
$\A$, and are all $\A$-equivalent. Moreover, these are maximum essential
quasi-multiplier extensions of $\A$.
\end{theorem}
\begin{proof}
The technique of the proof is parallel to that of \cite{BK} Theorem 3.2 and
\cite{B} Theorem 6.1.

Let $\pi: \A\to\B(\Hi)$ be a completely isometric homomorphism, and consider
the following canonical completely contractive homomorphisms:$$\A
\hookrightarrow\A^{**}\overset{\pi^{**}}{\to}\B(\Hi)^{**}\to\B(\Hi).$$ Let
$\hat{\pi}$ be the composition of the last two mappings. Then $\hat{\pi}$ is
completely isometric on $\cl{\hat{A}}$, and also $\hat{\pi}(\hat{a})=\pi(a),\;
\forall a\in\A$. Let $\tilde{\pi}:=\hat{\pi}|_{\QM^{**}(\A)}$.

Let us consider the canonical complete contractions:
$$\QM^{**}(\A)\overset{\tilde{\pi}}{\to}\QM^{\pi}(\A)\overset{\phi}{\to}
\QM^{\tilde{\psi}}(\A)\overset{\theta}{\to}\Q(\A),$$ where, $\tilde{\pi}$ and
$\phi$ are quasi-homomorphisms. Explicitly, $\phi=\rho\circ\sigma$, where
$\sigma$ is as in Theorem~\ref{situation} and $\rho$ is as in
Lemma~\ref{improvedqm}; $\theta(x)(a_1, a_2)=a_1\bullet v^*\bullet x\bullet
a_2$ for $x\in\QM^{\tilde{\psi}}(\A)$, $a_1, a_2\in\A$ with $v^*$ as in
Lemma~\ref{bk}.

To check that $\tilde{\pi}$ maps $\QM^{**}(\A)$ into $\QM^{\pi}(\A)$, take $z
\in\QM^{**}(\A)$, $a, b\in\A$. Then $\pi(a)\tilde{\pi}(z)\pi(b)=\hat{\pi}(\hat{
a}z\hat{b})\in\tilde{\pi}(\cl{\hat{A}})=\pi(\A)$. That $\tilde{\pi}$ is an
isometry follows from the facts that $\forall z\in\QM^{**}(\A),\;\;\|\tilde{\pi
}(z)\|\ge\|\pi(e_\alpha)\tilde{\pi}(z)\pi(e_\beta)\|=\|\hat{\pi}(\hat{e}_\alpha
z\hat{e}_\beta)\|=\|\hat{e}_\alpha z\hat{e}_\beta\|$, and that the identity of
$\A$ is a weak* limit point of $\{e_{\alpha}\}$, using the separate weak*
continuity of the product on $\A^{**}$. In fact, for each $\epsilon>0$ and each
$f\in Ball(\A^*)$, $\sup_{\alpha}\sup_{\beta}|(\hat{e}_{\alpha}z\hat{e}_{\beta}
)(f)|\ge|z(f)|-\epsilon$, so that $\sup_{\alpha}\sup_{\beta}\|\hat{e}_{
\alpha}z\hat{e}_{\beta}\|=\sup_{\alpha}\sup_{\beta}\sup_{f\in Ball(\A^*)}|(\hat
{e}_{\alpha}z\hat{e}_{\beta})|=\sup_{f\in Ball(\A^*)}\sup_{\alpha}\sup_{\beta}|
(\hat{e}_{\alpha}z\hat{e}_{\beta})|\ge\sup_{f\in Ball(\A^*)}|z(f)|-\epsilon=\|z
\|-\epsilon$. Since $\epsilon>0$ is arbitrary, $\sup_{\alpha}\sup_{\beta}\|\hat
{e}_{\alpha}z\hat{e}_{\beta}\|\ge\|z\|$. A similar calculation at the matrix
level shows that $\tilde{\pi}$ is a complete isometry.

Let $x\in\QM^{\tilde{\psi}}(\A)$, and write $\varphi_x:=\theta(x)$. Then for $a
, b, c\in\A$, $\varphi_x(a, b)=a\bullet v^*\bullet x\bullet b$, and $\varphi_x(
\cdot,a)\in CB(\A)$ and $\varphi_x(a, \cdot)\in CB(\A)$. Also $\varphi_x(ab,c)=
(ab)\bullet v^*\bullet x\bullet c=a\bullet v^*\bullet b\bullet v^*\bullet x
\bullet c=a(b\bullet v^*\bullet x\bullet c)=a\varphi_x(b, c)$. Similarly, $
\varphi_x(a,bc)=\varphi_x(a,b)c$. Hence $\theta(z)=\varphi_x\in\Q(\A)$. That $
\theta$ is one-to-one easily follows from \cite{BP}~Corollary~1.3. In fact, let
$\theta(x)=0$, i.e., $\varphi_x(a,b)=0,\;\;\forall a,b\in\A$. Then $a\bullet v^
*\bullet x\bullet b=0,\;\;\forall a,b\in\A$. By \cite{BP} Corollary 1.3, $v^*
\bullet x\bullet b=0$, so that $v\bullet v^*\bullet x\bullet b=0$, and hence $x
\bullet b=0,\;\;\forall b\in\A$ since $v\bullet v^*$ is the identity of $I_{11}
(\A)$. Thus again by \cite{BP} Corollary 1.3, $x=0$.

Let $\varphi\in\Q(\A)$, and let $F$ be a weak* accumulation point of
$\widehat{\varphi(e_\alpha, e_\alpha)}$ in $\A^{**}$. Clearly,
$\|F\|\le\|\varphi\|_{cb}$. For $a, b\in\A$, we have $\widehat{\varphi(a,
b)}=\lim_\alpha\widehat{\varphi(ae_\alpha, e_\alpha b)}=\lim_\alpha\hat{a}
\widehat{\varphi(e_\alpha, e_\alpha)}\hat{b}=\hat{a}F\hat{b}$. Hence, $F\in\QM^
{**}(\A)$ and $\theta\circ\phi\circ\tilde{\pi}(F)=\varphi$ and $\|\varphi\|_{cb
}\le\|F\|$, and thus $\theta\circ\phi\circ\tilde{\pi}$ is an onto isometry.
Here, to see that $\theta\circ\phi\circ\tilde{\pi}(F)=\varphi$, first note that
$\pi(a)\tilde{\pi}(F)\pi(b)=\pi(a\bullet\sigma(\tilde{\pi}(F))\bullet b)$ by
Theorem~\ref{situation}. But the left hand side is $\tilde{\pi}(\hat{a})\tilde{
\pi}(F)\tilde{\pi}(\hat{b})=\tilde{\pi}(\hat{a}F\hat{b})=\tilde{\pi}(\widehat{
\varphi(a,b)})=\pi(\varphi(a,b))$. On the other hand, $\pi(a\bullet\sigma(
\tilde{\pi}(F))\bullet b)=\pi(a\bullet v^*\bullet v\bullet\sigma(\tilde{\pi}(F)
)\bullet b)=\pi(a\bullet v^*\bullet\rho(\sigma(\tilde{\pi}(F)))\bullet b)=\pi(a
\bullet v^*\bullet\phi(\tilde{\pi}(F))\bullet b)=\pi(\theta(\phi(\tilde{\pi}(F)
))(a,b))$. Thus $\theta\circ\phi\circ\tilde{\pi}(F)=\varphi$ follows. Now since
we know that elements of $\Q(\A)$ are completely bounded, we can equip $\Q(\A)$
with a matrix norm as mentioned after the definition of $\Q(\A)$. Since the
operation of $\varphi$ is given by the multiplication by $F$, a similar
calculation works at the matrix level, and $\theta\circ\phi\circ\tilde{\pi}$
is a complete isometry. It is easy to check that $\phi(=\rho\circ\sigma)$ is
one-to-one by using the fact that $\pi(e_\alpha)\overset{SOT}{\to}1_{\Hi}$. In
fact, since $\rho$ is a complete isometry (Lemma \ref{improvedqm}), it suffices
to show that $\sigma$ is one-to-one. Let $\sigma(y)=0,\;\;y\in\QM^{\pi}(\A)$.
Then $\pi(e_{\alpha})y\pi(e_{\beta})=\pi(e_{\alpha}\bullet y\bullet e_{\beta})=
0$, thus $\pi(e_{\alpha})y\pi(e_{\beta})\xi=0,\;\forall \xi\in\Hi$. By taking
the limits $\alpha,\beta\to+\infty$, $y\xi=0,\;\forall \xi\in\Hi$, so that $y
=0$.

Thus, we have proved that $\tilde{\pi}$ is a complete isometry; $\phi$ and
$\theta$ are one-to-one complete contractions;
$\theta\circ\phi\circ\tilde{\pi}$ is an onto complete isometry. All these facts
force that each of $\tilde{\pi}$, $\phi$, and $\theta$ is an onto complete
isometry. Hence (I)-(III) are all $\A$-equivalent, and they are essential
quasi-multiplier extensions of $\A$, and $\Q(\A)$ is an operator space.

Finally, we prove that (III) is a maximum essential quasi-multiplier
extension. Let $(X, \pi)$ be any essential quasi-multiplier extension of $\A$.
Then by Theorem~\ref{situation} and Lemma~\ref{improvedqm}, there exists a
completely contractive quasi-homomorphism $\phi: X\to\QM^{\tilde{\psi}}(\A)$.
To see that $\phi\circ\pi=\tilde{\psi}$, take any $a,b,c\in\A$. By
Theorem~\ref{situation}, $\pi(a)\pi(b)\pi(c)=\pi(a\bullet\sigma(\pi(b))\bullet
c)$, so that $\pi(abc)=\pi(a\bullet\sigma(\pi(b))\bullet c)$, and thus $abc=a
\bullet\sigma(\pi(b))\bullet c$. But $abc=a\bullet v^*\bullet b\bullet v^*
\bullet c$ by Lemma~\ref{bk}~(3), so that $a\bullet\sigma(\pi(b))\bullet c=a
\bullet v^*\bullet b\bullet v^*\bullet c$. Since $a\in\A$ is arbitrary, by
\cite{BP}~Corollary~1.3, $\sigma(\pi(b))\bullet c=v^*\bullet b\bullet v^*
\bullet c$. Hence $\phi(\pi(b))\bullet c=\rho(\sigma(\pi(b)))\bullet c=v\bullet
\sigma(\pi(b))\bullet c=v\bullet v^*\bullet b\bullet v^*\bullet c=b\bullet v^*
\bullet c$ since $v\bullet v^*$ is the identity of $I_{11}(\A)$. Since $c\in\A$
is arbitrary, again by \cite{BP}~Corollary~1.3, $\phi(\pi(b))=b\bullet v^*=
\tilde{\psi}(b)$. Since $b\in\A$ is arbitrary, $\phi\circ\pi=\tilde{\psi}$.
Thus $(X,\pi)\le(\QM^{\tilde{\psi}}(\A), \tilde{\psi})$, so that $(\QM^{\tilde{
\psi}}(\A), \tilde{\psi})$ is a maximum essential quasi-multiplier extension.
\end{proof}
\begin{corollary}
In the C*-algebra case, our definition of the quasi-multipliers coincides with
the existing one (\cite{B}, also see \cite{Pe} \S 3.12) in the sense that
these are completely isometrically quasi-isomorphic.
\end{corollary}
\begin{corollary} If $\A$ is an operator algebra with a two-sided c.a.i., then
$\Q(\A)$ is an operator space, and it can be taken as a maximum essential
quasi-multiplier extension of $\A$.
\end{corollary}
\section{Quasi-multipliers and algebrizations of an operator space}
\label{section: algebrizations}
In this section we present the main theorem (Theorem~\ref{main}) of this paper.
We consider possible functors from the category of operator spaces
together with complete isometries into the category of operator algebras
together with completely isometric homomorphisms. Hence comes the name {\bf
algebrization}. But we do not have to express this ``algebrization functor''
explicitly, the name ``algebrization'' shows up only in the titles of this
paper and this section.

Theorem \ref{main} completely characterizes
an operator algebra
without any assumption on identities or approximate identities. The item (iii)
of the following theorem is regarded as the ``quasi'' version of the
Blecher-Effros-Zarikian Theorem (so-called ``$\tau$-trick'') (\cite{BEZ}
Theorem 4.6), that is, operator algebra products are characterized only in
terms of the matrix norm. To use the Haagerup tensor norm is essential. Also
a generalization of the Blecher-Ruan-Sinclair theorem is obtained as a simple
corollary.
\begin{theorem}\label{main} Let $X$ be a non-zero operator space with a
bilinear mapping $\varphi: X\times X\to X$, and let $I(\Sy_X)$ be as in Section
\ref{section: pre} and $1$ be its identity. We regard $X$ as a subspace of
$I(\Sy_X)$ by the \v{S}ilov embedding as explained in Section \ref{section:
pre}. Let
$$\begin{matrix} & \M_2\left(I(\Sy_X)\otimes_hI(\Sy_X)\right) & & \M_2(X)\\ &
\cup & & \cup\\\Gamma_{\varphi}: &
\left[\begin{matrix}X\otimes_h\C1 &
X\otimes_h X\\O & \C1\otimes_hX\end{matrix}\right] & \to &
\left[\begin{matrix}X & X \\ O & X\end{matrix}\right]\end{matrix}$$ be defined
by $$\Gamma_{\varphi}\left(\left[\begin{matrix}x_1\otimes1& x\otimes y\\0 &
1\otimes x_2\end{matrix}\right]\right):=\left[\begin{matrix}x_1 &
\varphi(x, y)\\0 & x_2\end{matrix}\right]$$ and their linear extensions
Then, the following are equivalent:
\begin{itemize}
\item[(i)] $(X, \varphi)$ is an abstract operator algebra (i.e., there is a
completely isometric homomorphism from $X$ into a concrete operator algebra,
hence, in particular, $\varphi$ is associative),
\item[(ii)] there exists a $z\in\QM(X)$ with $\|z\|\le1$ such that $\forall x,
y\in X$, $\varphi(x, y)=x\bullet z\bullet y$,
\item[(iii)] $\Gamma_{\varphi}$ is completely contractive.
\end{itemize}
Moreover, such a $z$ is unique.

When these conditions hold, we say that $\varphi$ is an {\bf operator algebra
product (OAP)} on $X$ and denote the set of all OAP's on $X$ by $OAP(X)$.
\end{theorem}
\begin{proof}\footnote{Historically, first I proved the equivalence of (i)
$\Leftrightarrow$(ii) separately, \cite{KP} Theorem~2.6.} Uniqueness of $z$
easily follows from \cite{BP}~Corollary~1.3. In fact, let $z_1,z_2\in\QM(X)$ be
such that $\|z_1\|\le1,\;\|z_2\|\le1,\;x\bullet z_1\bullet y=\varphi(x, y)=x
\bullet z_2\bullet y,\;\forall x,y\in X$. Then $x\bullet(z_1-z_2)\bullet y=0,\;
\forall x,y\in X$, so that $(z_1-z_2)\bullet y=0,\;\forall y\in X$ by
\cite{BP}~Corollary~1.3, and thus $(z_1-z_2)^*(z_1-z_2)\bullet y=0,\;\forall y
\in X$. Hence, again by \cite{BP}~Corollary~1.3, $(z_1-z_2)^*(z_1-z_2)=0$, and
$z_1=z_2$.

$\underline{\text{(ii)}\implies\text{(i)}:}$ This direction follows
from the same way as Remark 2 on page 194 of \cite{BRS}. Define $\rho:X\to
I(\Sy_X)$ by $$\rho(x):=\left[\begin{matrix}x\bullet z &
x\bullet\sqrt{1_{22}-z\bullet z^*}\\0 & 0\end{matrix}\right],\;\;\forall
x\in X.$$ Since
\begin{align*} & \rho(x)=\left[\begin{matrix}0 & x\\0 &
0\end{matrix}\right]\odot\left[\begin{matrix}0 & 0\\z &
\sqrt{1_{22}-z\bullet z^*}\end{matrix}\right]\;\;\text{and}\\
& \left[\begin{matrix}0 & 0\\z &
\sqrt{1_{22}-z\bullet z^*}\end{matrix}\right]\odot\left[\begin{matrix}0
& z^*\\0 & \sqrt{1_{22}-z\bullet z^*}\end{matrix}\right]=\left[\begin{matrix}0
& 0\\0 & 1_{22}\end{matrix}\right],\end{align*} it follows that $\rho$ is a
completely isometric homomorphism.

$\underline{\text{(i)}\implies\text{(iii)}:}$ We may assume that $(X,
\varphi)\subset(\B(\Hi), \cdot)$ as a subalgebra of operators. By the
construction of $I(\Sy_X)$ in Section \ref{section: pre},
$1_{11}=\left[\begin{matrix}1_{11} & 0\\0 &
0\end{matrix}\right]=\left[\begin{matrix}1_{\Hi} & 0\\0 &
0\end{matrix}\right]$, $1_{22}=\left[\begin{matrix}0 & 0\\0 &
1_{22}\end{matrix}\right]=\left[\begin{matrix}0 & 0\\0 &
1_{\Hi}\end{matrix}\right]$. We write $1_{21}:=\left[\begin{matrix}0 &
0\\1_{\Hi} & 0\end{matrix}\right]\in\M_2(\B(\Hi))$ and define a linear
mapping\\$\phi: I(\Sy_X)\otimes_hI(\Sy_X)\to\M_2(\B(\Hi))$ by
$\phi(\eta\otimes\zeta):=\eta\cdot1_{21}\cdot\zeta\;\;\forall \eta, \zeta\in
I(\Sy_X)$, and their linear extensions, where the products between
elements of $\M_2(\B(\Hi))$ are induced from the original products $\cdot$ of
$\B(\Hi)$, and still denoted by ``\;$\cdot$\;''. Then, obviously $\phi$ is
completely contractive. Let $$\xi=\left[\begin{matrix}x\otimes1 &
\sum_{i=1}^nx_i\otimes y_i\\0 & 1\otimes y\end{matrix}\right]
\in\left[\begin{matrix}X\otimes_h\C1 & X\otimes_hX\\O & \C1\otimes_hX
\end{matrix}\right],$$ then
\begin{align*}\|\Gamma_\varphi(\xi)\| & =\left\|\Gamma_\varphi\left(\left[
\begin{matrix}x\otimes1 &
\sum_{i=1}^nx_i\otimes y_i\\0 & 1\otimes y\end{matrix}\right]\right)\right\|
=\left\|\left[\begin{matrix}x & \sum_{i=1}^n\varphi(x_i, y_i)\\0 & y
\end{matrix}\right]\right\|\\
& =\left\|\left[\begin{matrix}x\cdot1_{21}\cdot1 &
\sum_{i=1}^nx_i\cdot1_{21}\cdot y_i\\0\cdot1_{21}\cdot0 & 1\cdot1_{21}\cdot
y\end{matrix}\right]\right\|=\left\|\left[\begin{matrix}
\phi(x\otimes1) & \phi(\sum_{i=1}^nx_i\otimes y_i)\\\phi(0\otimes0) &
\phi(1\otimes y)\end{matrix}\right]\right\|\\
& =\left\|\phi_2\left(\left[\begin{matrix}x\otimes1 & \sum_{i=1}^nx_i\otimes
y_i\\0 & 1\otimes y\end{matrix}\right]\right)\right\|=\|\phi_2(\xi)\|,
\end{align*}where $\phi_2: \M_2(I(\Sy_X))\to\M_2(\M_2(\B(\Hi)));\;
\phi_2((x_{i, j})):=(\phi(x_{i, j}))$. Since $\phi$ is completely contractive,
we obtain that $\|\Gamma_\varphi(\xi)\|\le\|\xi\|$, so that $\Gamma_\varphi$ is
contractive. A similar calculation at the matrix level shows that
$\Gamma_\varphi$ is completely contractive.

Now, we show the hardest direction.

$\underline{\text{(iii)}\implies\text{(ii)}:}$ Let $\widetilde{I(\Sy_X)}$ be a
copy of $I(\Sy_X)$ which shares $\C1$ with $I(\Sy_X)$, and let $\;\widetilde{}
\;: I(\Sy_X)\to\widetilde{I(\Sy_X)}$ be the canonical mapping. Note that
$\tilde{1}=1$. $I(\Sy_X)*_1\widetilde{I(\Sy_X)}$ be the completion of the free
product of $I(\Sy_X)$ and $\widetilde{I(\Sy_X)}$ amalgamated over $\C1$, which
is a C*-algebra. We embed $I(\Sy_X)\otimes_hI(\Sy_X)$ into
$I(\Sy_X)*_1\widetilde{I(\Sy_X)}$ by the complete isometry $\gamma$ defined by
setting $\gamma(\xi\otimes\eta):=\xi*\tilde{\eta}$ and extending linearly. The
reader unfamiliar with this embedding is recommended to refer to \cite{Pi0},
\cite{Pi}, and also \cite{P} Chapter 17. The important properties of this free
product employed in the proof below are that $I(\Sy_X)*_1\widetilde{I(\Sy_X)}$
contains both $I(\Sy_X)$ and $\widetilde{I(\Sy_X)}$ as C*-subalgebras, and
that these three C*-algebras share the common identity
$1=\tilde{1}=1*\tilde{1}$. Let us define
$$S:=\left[\begin{matrix}\C1_{11} & \gamma(X\otimes_h\C1) & O &
\gamma(X\otimes_hX)\\\gamma(X\otimes_h\C1)^* & \C1_{22} & O & O\\O & O &
\C\tilde{1}_{11} & \gamma(\C1\otimes_hX)\\\gamma(X\otimes_hX)^* & O &
\gamma(\C1\otimes_hX)^* & \C\tilde{1}_{22}\end{matrix}\right],$$which is a
subset of $\M_4(I(\Sy_X)*_1\widetilde{I(\Sy_X)})$. Let $C^*(S)$ be the
C*-algebra generated by $S$ in $\M_4(I(\Sy_X)*_1\widetilde{I(\Sy_X)})$. The
elements $[\xi_{i, j}]_{1\le i, j\le4}\in C^*(S)$ of the form $\xi_{i,
j}=\xi_{i, j}^1*\cdots*\xi_{i, j}^{n_{i, j}},\;n_{i, j}\in\N$ with$$\xi_{i,
j}^1\in\begin{cases}X, & \text{if $i=1$;}\\X^*, & \text{if $i=2$;}\\
\widetilde{X}, & \text{if $i=3$;}\\\widetilde{X}^*, & \text{if $i=4$,}
\end{cases}\quad\quad\text{and}\quad\quad\xi_{i, j}^{n_{i,
j}}\in\begin{cases}X^*, & \text{if $j=1$;}\\X, & \text{if $j=2$;}\\
\widetilde{X}^*, & \text{if $j=3$;}\\\widetilde{X}, & \text{if $j=4$}
\end{cases}$$span a dense subset of $C^*(S)$. Hence $diag\{1_{11}, 1_{22},
\tilde{1}_{11}, \tilde{1}_{22}\}$ is the identity of $C^*(S)$, so that $S$ is
an operator system in $C^*(S)$. Define a linear mapping
$$\Phi: S\to\left[\begin{matrix}\C1_{11} & X & O & X\\X^* & \C1_{22} &
O & O\\O & O & \C1_{11} & X\\X^* & O & X^* & \C1_{22}\end{matrix}\right]
\subset\M_2(I(\Sy_X))$$by\begin{align*}\Phi\left(\left[\begin{matrix}
\lambda1_{11} & \gamma(x_1\otimes1) & 0 & \gamma(x_5\otimes x_6)\\
\gamma(x_2\otimes1)^* & \mu1_{22} & 0 & 0\\0 & 0 & \alpha\tilde{1}_{11} &
\gamma(1\otimes x_3)\\\gamma(x_7\otimes x_8)^* & 0 & \gamma(1\otimes x_4)^* &
\beta\tilde{1}_{22}\end{matrix}\right]\right)\\:=\left[\begin{matrix}
\lambda1_{11} & x_1 & 0 & \varphi(x_5, x_6)\\x_2^* & \mu1_{22} & 0 & 0\\0 & 0 &
\alpha1_{11} & x_3\\\varphi(x_7, x_8)^* & 0 & x_4^* & \beta1_{22}\end{matrix}
\right]\end{align*}and their linear extensions. By the canonical shuffle, and
the fact that $\gamma$ is a complete isometry, and our assumption that
$\Gamma_{\varphi}$ is completely contractive, we know that $\Phi$ is completely
positive. We extend $\Phi$ to a linear mapping $\Phi'$ from $$S':=span\{S\cup
diag\{\C1, \C1, \C1, \C1\}\}$$ onto the same range such that
$$\Phi'(diag\{\lambda1, \mu1, \alpha1, \beta1\})=diag\{\lambda1_{11},
\mu1_{22}, \alpha1_{11}, \beta1_{22}\}.$$ This uniquely well defines $\Phi'$.
In particular,$$\text{Ker}\Phi'=\{\C1_{22}, \C1_{11}, \C\tilde{1}_{22},
\C\tilde{1}_{11}\}.$$ To see that $\Phi'$ is completely positive, simply
observe that for $\xi\in S'$,$$\Phi'(\xi)=\Phi(diag\{1_{11}, 1_{22},
\tilde{1}_{11}, \tilde{1}_{22}\}\;\xi\;diag\{1_{11}, 1_{22}, \tilde{1}_{11},
\tilde{1}_{22}\}).$$ Since $S'$ is an operator system containing the identity
of $\M_4(I(\Sy_X)*_1\widetilde{I(\Sy_X)})$, we can extend $\Phi'$ to a
completely positive map$$\widetilde{\Phi}:
\M_4(I(\Sy_X)*_1\widetilde{I(\Sy_X)})\to\left[\begin{matrix}
I_{11}(X) & I(X) & I_{11}(X) & I(X) \\
I(X)^* & I_{22}(X) & I(X)^* & I_{22}(X) \\
I_{11}(X) & I(X) & I_{11}(X) & I(X) \\
I(X)^* & I_{22}(X) & I(X)^* & I_{22}(X)\end{matrix}\right]$$ by the
injectivity of the right hand side.
Since $\widetilde{\Phi}$ ``fixes'' each scalar element on the diagonal,
$\widetilde{\Phi}$ is factored to $[\phi_{i, j}]_{1\le i, j\le 4}$ by a
common argument. $\widetilde{\Phi}$ also ``fixes'' the C*-subalgebra
$$\left[\begin{matrix}I_{11}(X) & I(X) & O & O\\I(X)^* & I_{22}(X) & O & O\\O &
O & \widetilde{I_{11}(X)} & \widetilde{I(X)}\\
O & O & \widetilde{I(X)^*} & \widetilde{I_{22}(X)}\end{matrix}\right]$$
by the rigidity. Hence, $\widetilde{\Phi}$ is a ``module map'' over it in the
sense of \cite{BP} Lemma~1.6.
Let $x, y\in X$, then
\begin{align*}
& {\left[\begin{matrix}
0 & 0 & 0 & \varphi(x, y)\\
0 & 0 & 0 & 0\\
0 & 0 & 0 & 0\\
0 & 0 & 0 & 0\end{matrix}\right]}
{=\widetilde{\Phi}\left(\left[\begin{matrix}
0 & 0 & 0 & x*\tilde{y}\\
0 & 0 & 0 & 0\\
0 & 0 & 0 & 0\\
0 & 0 & 0 & 0\end{matrix}\right]\right)}\\
= & {\widetilde{\Phi}\left(\left[\begin{matrix}
0 & x & 0 & 0\\
0 & 0 & 0 & 0\\
0 & 0 & 0 & 0\\
0 & 0 & 0 & 0\end{matrix}\right]\left[\begin{matrix}
0 & 0 & 0 & 0\\
0 & 0 & 1 & 0\\
0 & 0 & 0 & 0\\
0 & 0 & 0 & 0\end{matrix}\right]\left[\begin{matrix}
0 & 0 & 0 & 0\\
0 & 0 & 0 & 0\\
0 & 0 & 0 & \tilde{y}\\
0 & 0 & 0 & 0\end{matrix}\right]\right)}\\
= & {\left[\begin{matrix}
0 & x & 0 & 0\\
0 & 0 & 0 & 0\\
0 & 0 & 0 & 0\\
0 & 0 & 0 & 0\end{matrix}\right]\widetilde{\Phi}\left(\left[\begin{matrix}
0 & 0 & 0 & 0\\
0 & 0 & 1 & 0\\
0 & 0 & 0 & 0\\
0 & 0 & 0 & 0\end{matrix}\right]\right)\left[\begin{matrix}
0 & 0 & 0 & 0\\
0 & 0 & 0 & 0\\
0 & 0 & 0 & y\\
0 & 0 & 0 & 0\end{matrix}\right]}\\
= & {\left[\begin{matrix}
0 & x & 0 & 0\\
0 & 0 & 0 & 0\\
0 & 0 & 0 & 0\\
0 & 0 & 0 & 0\end{matrix}\right]\left[\begin{matrix}
0 & 0 & 0 & 0\\
0 & 0 & z & 0\\
0 & 0 & 0 & 0\\
0 & 0 & 0 & 0\end{matrix}\right]\left[\begin{matrix}
0 & 0 & 0 & 0\\
0 & 0 & 0 & 0\\
0 & 0 & 0 & y\\
0 & 0 & 0 & 0\end{matrix}\right]=\left[\begin{matrix}
0 & 0 & 0 & x\bullet z\bullet y\\
0 & 0 & 0 & 0\\
0 & 0 & 0 & 0\\
0 & 0 & 0 & 0\end{matrix}\right],}
\end{align*}
where $z:=\phi_{2, 3}(1)$.
\end{proof}
The new point in the following corollary is that we do not assume that $(X,
\varphi)$ has a ``two-sided'' c.a.i.. This also could follow from \cite{B0}.
\begin{corollary}(A generalization of the BRS Theorem \cite{BRS}, \cite{R})
Let $X$ be a non-zero operator space, and $\varphi$ be a completely contractive
bilinear mapping on $X$. If there exists $\{e_\alpha\}$ and $\{f_\beta\}$ with
$\|e_\alpha\|\le1,\;\|f_\beta\|\le1$ such that $\lim_\alpha\varphi(x, e_\alpha)
=\lim_\beta\varphi(f_\beta, x)=x,\;\forall x\in X$,\footnote{In this case, $(X,
\varphi)$ has a two-sided contractive approximate identity.} then $(X, \varphi)
$ is an operator algebra.
\end{corollary}
\begin{proof} For any element $\xi=\left[\begin{matrix}x\otimes1 &
\sum_{i=1}^nx_i\otimes y_i\\0 & 1\otimes y\end{matrix}\right]
\in\left[\begin{matrix}X\otimes_h\C1 & X\otimes_hX\\O & \C1\otimes_hX
\end{matrix}\right]$,
\begin{align*}& 
\Gamma_\varphi\left(\left[\begin{matrix}x\otimes1 & \sum_{i=1}^nx_i\otimes y_i
\\0 & 1\otimes y\end{matrix}\right]\right)=\left[\begin{matrix}x &
\sum_{i=1}^n\varphi(x_i, y_i)\\0 & y\end{matrix}\right]\\
= & \lim_\alpha\lim_\beta\left[\begin{matrix}
\varphi(x, e_\alpha) & \sum_{i=1}^n\varphi(x_i, y_i)\\0 & \varphi(f_\beta, y)
\end{matrix}\right]=\lim_\alpha\lim_\beta\tilde{\varphi}_2\left(\left[
\begin{matrix}x\otimes e_\alpha & \sum_{i=1}^nx_i\otimes y_i\\0 &
f_\beta\otimes y\end{matrix}\right]\right)\\
= & \lim_\alpha\lim_\beta\tilde{\varphi}_2\left(\left[
\begin{matrix}1\otimes1 & 0\\0 & f_\beta\otimes1\end{matrix}\right]
\left[\begin{matrix}x\otimes1 & \sum_{i=1}^nx_i\otimes y_i\\0 & 1\otimes y
\end{matrix}\right]\left[\begin{matrix}1\otimes e_\alpha & 0\\0 & 1\otimes1
\end{matrix}\right]\right),
\end{align*} where $\tilde{\varphi}$ is the linear mapping $X\otimes_hX\to X$
corresponding to the bilinear mapping $\varphi: X\times X\to X$, and
$\tilde{\varphi}_2: \M_2(X\otimes_hX)\to\M_2(X);\;\tilde{\varphi}_2((\zeta_{i,
j})):=(\varphi(\zeta_{i, j}))$. Also we abused notation, and regarded $I(\Sy_X)
\otimes_hI(\Sy_X)$ as a subset of $I(\Sy_X)*_1I(\Sy_X)$ by the canonical
complete isometry as in the proof of Theorem~\ref{main} (iii)$\implies$(ii),
and the product between elements of $I(\Sy_X)\otimes_hI(\Sy_X)$ is taken in
$I(\Sy_X)*_1I(\Sy_X)$. Hence,
$\|\Gamma_\varphi(\xi)\|\le\|\tilde{\varphi}_2\|\|\xi\|
\le\|\tilde{\varphi}\|_{cb}\|\xi\|\le\|\xi\|$, so that $\Gamma_\varphi$ is
contractive. Similarly, $\Gamma_\varphi$ is completely contractive. Hence, by
\ref{main}, $(X, \varphi)$ is an operator algebra.
\end{proof}

We close this paper by providing an equivalent condition for
$\|z\|=\|\varphi\|_{cb}$ to hold. This is a simple corollary of Theorem
\ref{main}.
\begin{corollary} Let $X$ be an operator space, $\varphi$ be a bilinear mapping
on $X$. Under the same notation as Theorem \ref{main}, let
$$\begin{matrix} & \M_2\left(I(\Sy_X)\otimes_hI(\Sy_X)\right) & & \M_2(X)\\
 & \cup & & \cup \\
\tilde{\Gamma}_{\varphi}: & \left[\begin{matrix}X\otimes_h\C1 & X\otimes_h X\\
O & \C1\otimes_hX\end{matrix}\right] & \to & \left[\begin{matrix}X & X \\ O &
X\end{matrix}\right]\end{matrix}$$ be defined by
$$\tilde{\Gamma}_{\varphi}\left(\left[\begin{matrix} x_1\otimes 1 & x\otimes
y\\0 & 1\otimes
x_2\end{matrix}\right]\right):=\left[\begin{matrix}\|\varphi\|_{cb}x_1 &
\varphi(x, y)\\0 & \|\varphi\|_{cb}x_2\end{matrix}\right]$$
and their linear extension. Then the following are equivalent:
\begin{itemize}
\item[(i)] there exists a $z\in\QM(X)$ with $\|z\|=\|\varphi\|_{cb}$ such
that $\forall x, y\in X$ $\varphi(x, y)=x\bullet z\bullet y$,
\item[(ii)] $\|\tilde{\Gamma}_{\varphi}\|_{cb}=\|\varphi\|_{cb}$.
\end{itemize}
Moreover, such a $z$ is unique.
\end{corollary}
\begin{proof}Uniqueness of $z$ follows in the same way as Theorem \ref{main}.

$\underline{\text{(i)}\implies\text{(ii)}:}$ Note that by the definition of
$\tilde{\Gamma}_{\varphi}$, clearly
$\|\tilde{\Gamma}_{\varphi}\|_{cb}\ge\|\varphi\|_{cb}$. To show that
$\|\Gamma_{\varphi}\|\le\|\varphi\|_{cb}$, we may assume that
$\|\varphi\|_{cb}=1$. Then $\|\Gamma_{\varphi}\|_{cb}\le1$ follows from Theorem
\ref{main} (ii)$\implies$(iii).

$\underline{\text{(ii)}\implies\text{(i)}:}$ By scaling, it is enough to
consider the case $\|\varphi\|_{cb}=1$. By Theorem
\ref{main} (iii)$\implies$(ii), there exists $z\in\QM(X)$ with $\|z\|\le 1$
such that $\forall x, y\in X$ $\varphi(x, y)=x\bullet z\bullet y$. Thus
$\|z\|\le\|\varphi\|_{cb}$. But clearly, $\|z\|\ge\|\varphi\|_{cb}$.
\end{proof}
Note that the difference between $\Gamma_{\varphi}$ and
$\tilde{\Gamma}_{\varphi}$ is ``more than scaling''.

\vspace{4mm}
Acknowledgments: This work was carried out while I was a graduate student of
Professor Vern I. Paulsen at University of Houston, and the work is part of my
Ph.D. thesis, August 2003 (\cite{K}). I am grateful to him for reading the
manuscript carefully with criticism, and for constant encouragement. For
the result of Theorem \ref{main}, I was awarded from the Sigma Xi, April 2003.


\vspace{4 mm}

\begin{thebibliography}{99}
\bibitem{B0}  D. P. Blecher,
{\em The Shilov boundary of an operator space and the characterization
  theorems}, Journal of Functional Analysis {\bf 182} (2001), 280-343.
\bibitem{B}  D. P. Blecher,
{\em One-sided ideals and approximate identities in operator algebras}, to
appear in Journal of the Australian Mathematical Society.
\bibitem{BEZ} D. P. Blecher, E. G. Effros and V. Zarikian,
{\em One-sided M-ideals and multipliers in operator spaces, I}, Pacific Journal
of Mathematics {\bf 206 (2)} (2002), 287-319.
\bibitem{BK} D. P. Blecher and M. Kaneda, {\em The ideal envelope of an
operator algebra}, to appear in Proceedings of the American Mathematical
Society.
\bibitem{BP} D. P. Blecher and V. I. Paulsen, {\em Multipliers of operator
spaces, and the injective envelope}, Pacific Journal of Mathematics {\bf 200}
(2001), 1-17.
\bibitem{BRS} D. P. Blecher, Z.-J. Ruan and A. M. Sinclair, {\em A
characterization of operator algebras}, Journal of Functional Analysis {\bf 89}
(1990), 188-201.
\bibitem{Br} L. G. Brown, {\em Close hereditary C*-algebras and the structure
of quasi-multipliers}, Mathematical Sciences Research Institute preprint, 1985.
\bibitem{CE} M.-D. Choi and E. G. Effros, {\em Injectivity and Operator
Spaces}, Journal of Functional Analysis {\bf 24} (1977), 156-209.
\bibitem{CS}  E. Christensen and A. M. Sinclair, {\em Representations of
completely bounded multilinear operators}, Journal of Functional Analysis {\bf
72} (1987), 151-181.
\bibitem{ER} E. G. Effros and Z.-J. Ruan, {\em Operator Spaces}, Oxford, 2000.
\bibitem{K} M. Kaneda {\em Multipliers and Algebrizations of Operator Spaces},
Ph.D. Thesis, University of Houston, August 2003.
\bibitem{KP} M. Kaneda and V. I. Paulsen,
{\em Quasimultipliers of operator spaces}, preprint 2003.
\bibitem{P} V. I. Paulsen,
{\em Completely bounded maps and operator algebras}, Cambridge University
Press, 2002.
\bibitem{Pe} G. K. Pedersen,
{\em C*-Algebras and their Automorphism Groups}, Academic Press, 1979.
\bibitem{Pi0} G. Pisier, {\em A simple proof of a theorem of Kirchberg and
related results on $C\sp *$-norms}, Journal of Operator Theory {\bf 35 (2)}
(1996), 317-335.
\bibitem{Pi} G. Pisier, {\em An Introduction to the Theory of Operator Spaces},
Cambridge University Press, 2003.
\bibitem{R} Z.-J. Ruan, {\em A characterization of non-unital operator
algebras}, Proceedings of the American Mathematical Society {\bf 121} (1994),
193-198.
\end{thebibliography}
\end{document}